\documentclass[10pt]{amsart}

\usepackage{amssymb, amsthm}
\newtheorem{thm}{Theorem}[section]
\newtheorem{lem}[thm]{Lemma}

\newtheorem{prop}[thm]{Proposition}

\theoremstyle{definition}
\newtheorem{defn}[thm]{Definition}
\newtheorem{ex}[thm]{Example}

\theoremstyle{remark}
\newtheorem{rem}[thm]{Remark}

\newcommand\CM{Cohen--Macaulay }

\def\NE{non--equidimensional }
\def\bs{^{\bigstar}}
\def\dot{^{\bullet}}
\def\HR{\hat R}

\begin{document}
\large

\title{Tight closure in non--equidimensional rings}

\author{Anurag K. Singh}

\address{Department Of Mathematics, University Of Michigan, East Hall, \\
525 East University Avenue, Ann Arbor, MI 48109-1109}

\maketitle

\section{Introduction}

Throughout our discussion, all rings are commutative, Noetherian and have an 
identity element. The notion of the {\it tight closure} of an ideal was 
developed by M.~Hochster and C.~Huneke in \cite{HHjams} and has yielded many 
elegant and powerful results in commutative algebra. The theory leads to the 
notion of {\it F--rational} rings, defined by R.~Fedder and K.-i.~Watanabe as 
rings in which parameter ideals are tightly closed, see \cite{FW}. Over a field
of characteristic zero, rings of F--rational type are now known to be precisely 
those having rational singularities by the work of K.~E.~Smith and N.~Hara, see 
\cite{Sminv, Hara}.

\medskip

We begin by recalling a theorem of Hochster and Huneke which states that a 
local ring $(R,m)$ which is 
a homomorphic image of a Cohen--Macaulay ring is F--rational if and only if it 
is equidimensional and has a system of parameters which generates a tightly 
closed ideal, \cite[Theorem 4.2 (d)]{HHbasec}. This leads to the question 
of whether a local ring in which a single 
system of parameters generates a tightly closed ideal must be equidimensional 
(and hence F--rational), \#19 of Hochster's \lq\lq Twenty Questions\rq\rq \ in 
\cite{Ho}. Rephrased, can a \NE ring have a system of parameters which
generates a tightly closed ideal -- we show it cannot for some classes of \NE 
rings. 

\medskip

A key point is that in equidimensional rings, tight closure has the so--called 
\lq\lq colon capturing\rq\rq \ property. This property does not hold in 
non--equidimensional rings. A study of these issues leads to a new 
closure operation, that we call {\it NE closure}. This closure does possess the 
colon capturing property even in non--equidimensional rings, and agrees with 
tight closure when the ring is equidimensional. We shall show that an excellent 
local ring $R$ is F--rational if and only if it has a system of parameters 
which generates an NE--closed ideal.

\section{Notation and terminology}

Let $R$ be a Noetherian ring of characteristic $p > 0$. We shall always use
the letter $e$ to denote a variable nonnegative integer, and $q$ to denote the
$e\,$th power of $p$, i.e., $q=p^e$. For an ideal 
$I=(x_1, \dots, x_n) \subseteq R$, we let $I^{[q]}=(x_1^q, \dots, x_n^q)$. We
shall denote by $R^{\text{o}}$ the complement of the union of the minimal 
primes of $R$.
For an ideal $I \subseteq R$ and an element $x$ of $R$, we say that 
$x \in I^*$, the {\it tight closure }\/ of $I$, if there exists 
$c \in R^{\text{o}}$ such 
that $cx^q \in I^{[q]}$ for all $q=p^e \gg 0$. If $I=I^*$, we say $I$ is 
{\it tightly closed}.

\medskip

An ideal $I \subseteq R$ is said to be a {\it parameter ideal}\/  if 
$I=(x_1,\dots,x_n)$ such that the images of $x_1,\dots,x_n$ form a system of
parameters in $R_P$, for every prime $P$ containing $I$.
The ring $R$ is said to be {\it F--rational}\/ if every parameter ideal of $R$ 
is tightly closed. 

\medskip

We next recall some well known results.

\begin{thm}

\item $(1)$\quad An F--rational ring $R$ is normal. If in addition $R$ is 
assumed to be the homomorphic image of a \CM ring, then $R$ is Cohen--Macaulay.

\item $(2)$\quad A local ring $(R,m)$ which is the homomorphic image of a 
\CM ring is F--rational if and only if it is equidimensional and the ideal
generated by one system of parameters is tightly closed. 

\item $(3)$\quad Let $P_1,\dots,P_n$ be the minimal primes of the ring $R$. 
Then for an 
ideal $I \subseteq R$ and $x \in R$, we have $x \in I^*$ if and only if for 
$1 \le i \le n$, its image $\overline x$ is in $(IR/P_i)^*$, the tight 
closure here being computed in the domain $R/P_i$.
\label{longlist}
\end{thm}

\begin{proof}
(1) and (2) are part of \cite[Theorem 4.2]{HHbasec} and (3) is observed as
\cite[Proposition 6.25 (a)]{HHjams}.
\end{proof}

\section{Main results}

\begin{lem} 
Let $P_1,\dots, P_n$ be the minimal primes of the ring 
$R$. Then for a tightly closed ideal $I$, we have $I =\bigcap_{i=1}^n (I+P_i)$.
\label{intersect}
\end{lem}

\begin{proof}  
That $I \subseteq \bigcap_{i=1}^n (I+P_i)$, is trivial. For the other
containment note that if 
$x \in \bigcap_{i=1}^n (I+P_i)$, then $\overline x \in (IR/P_i)^*$ for 
$1 \le i \le n$. Now by Theorem \ref{longlist} (3), we get that 
$x \in I^*=I$.
\end{proof}

The following theorem, although it has some rather strong hypotheses, does
show that no ideal generated by a system of parameters is tightly closed in the
\NE rings 
$$
R=K[[X_1,\dots, X_n,Y_1,\dots, Y_m]]/(X_1, \dots, X_n)\cap
(Y_1, \dots, Y_m)
$$ 
where $m,n \ge 1$ and $m \neq n$. 

\begin{thm} 
Let $(R,m)$ be a \NE local ring, with the minimal primes
partitioned into the sets $\{P_i\}$ and $\{Q_j\}$, such that 
$\dim R = \dim R/P_i > \dim R/Q_j$, for all $i$ and $j$. Let $P$ and $Q$ be the
intersections, $P=\bigcap P_i$ and $Q=\bigcap Q_j$. If $I \subseteq P+Q$ is an 
ideal of $R$ which is generated by a 
system of parameters, then $I$ cannot be tightly closed. In particular, if 
$P+Q=m$, then no ideal of $R$ generated by a system of parameters is tightly 
closed.
\label{contain}
\end{thm}

\begin{proof}  
Suppose not, let $I=(p_1+q_1, \dots, p_n+q_n)$ be a tightly closed 
ideal of $R$, where $p_1+q_1, \dots, p_n+q_n$ is a system of parameters with 
$p_i \in P$ and $q_i \in Q$. Note that we have
$$
(I+P)\cap (I+Q) \subseteq (\cap (I+P_i))\cap (\cap(I+Q_j))=I
$$
using Lemma \ref{intersect}. Consequently $(I+P)\cap(I+Q)
\subseteq I$, and so $p_i,q_i \in I$. In particular, 
$p_i = r_1(p_1+q_1)+ \dots +r_n(p_n+q_n)$. We first note that if $r_i \notin m$ 
then $q_i \in (p_1+q_1, \dots, p_{i-1}+q_{i-1},p_i,
                          p_{i+1}+q_{i+1}, \dots, p_n+q_n)$, but then $p_i\in P$
is a parameter, which is impossible since $\dim R = \dim R/P$.
Hence $r_i \in m$, and so $1-r_i$ is a unit. From this we may conclude that
$p_i \in (p_1+q_1, \dots, p_{i-1}+q_{i-1},q_i,
                          p_{i+1}+q_{i+1}, \dots, p_n+q_n)$, and so the ideal
$I$ may be written as
$I=(p_1+q_1, \dots, p_{i-1}+q_{i-1},q_i,
                          p_{i+1}+q_{i+1}, \dots, p_n+q_n)$. 

Proceeding this way, we see that $I=(q_1, \dots, q_n)$, i.e., $I\subseteq
Q$. But then each $Q_j=m$, a contradiction.
\end{proof}
 
\begin{rem}                                             
We would next like to discuss briefly the case where the \NE local ring $(R,m)$ 
is 
of the form  $R=S/(P \cap Q)$ where $S$ is a regular local ring with primes $P$
and $Q$ of different height. Then $R$ has minimal primes $\overline P$ and 
$\overline Q$ where, without loss of generality, 
$\dim R = \dim R/\overline P > \dim R/\overline Q$. If $I$ is an ideal of $S$
whose image $\overline I$ in $R$ is a tightly closed ideal, we see that 
$I + (P \cap Q) = (I+P)\cap (I+Q)$ by Lemma \ref{intersect},
and so it would certainly be enough to show that this cannot hold when 
$\overline I$ is generated by a system of parameters for $R$. One can indeed 
prove this in the case $S/P$ is Cohen--Macaulay, and $S/Q$ is a discrete 
valuation ring, which is Theorem \ref{CMDVR} below. However if we drop the 
hypothesis that $S/P$ be Cohen--Macaulay, this is no longer true: see 
Example \ref{badintersect}.
\end{rem}

\begin{lem} 
If $I$, $P$, and $Q$ are ideals of $S$, satisfying the condition that
$I+(P \cap Q)=(I+P) \cap (I+Q)$, 
then $I \cap (P+Q)=(I \cap P)+(I \cap Q)$.
\label{distribute}
\end{lem}

\begin{proof}  
Let $i=p+q \in I \cap (P+Q)$, where $p \in P$ and $q \in Q$. Then
$i-p=q \in (I+P) \cap (I+Q)=I+(P \cap Q)$ and so $i-p=q=\tilde i +r$ where
$\tilde i \in I$ and $r \in P \cap Q$. Finally note that
$i=(i - \tilde i)+\tilde i \in (I \cap P)+(I \cap Q)$, since 
$i - \tilde i = p+r \in I \cap P$ and $\tilde i=q-r \in I \cap Q$.
\end{proof}

\begin{lem} 
Let $M$ be an $S$--module and $N$, a submodule. If 
$x_1, \dots, x_n$ are elements of $S$ which form a regular sequence on $M/N$,
then 
$$
(x_1, \dots, x_n)M \cap N = (x_1, \dots, x_n)N.
$$

In particular, if $I$ and $J$ are ideals of $S$ and $I$ is generated by elements
which form a regular sequence on $S/J$, then $I \cap J = IJ$.
\label{regseq}
\end{lem}

\begin{proof}  We shall proceed by induction on $n$, the number of elements. If $n=1$,
the result is simple. For the inductive step, note that if we have 
$u=x_1m_1+ \dots +x_km_k \in (x_1, \dots, x_k)M \cap N$ with $m_i \in M$, then 
since $x_k$ is not a zero divisor on the module $M/((x_1, \dots, x_{k-1})M+N)$, 
we get that $m_k \in (x_1, \dots, x_{k-1})M+N$. Consequently, 
$u \in (x_1, \dots, x_{k-1})M +x_kN$. 
\end{proof}

\begin{thm}
Let $S=K[[X_1, \dots, X_n, Y]]$ with ideal $Q=(X_1,\dots, X_n)S$, and 
ideal $P$ satisfying the condition that $S/P$ is Cohen--Macaulay. Then if 
$R=S/(P \cap Q)$ is a \NE ring, no ideal of $R$ generated by a system of 
parameters can be tightly closed.
\label{CMDVR}
\end{thm}

\begin{proof} 
Let $I$ be an ideal of $S$ generated by elements which map to a system of
parameters in $R$. If the image of $I$ is a tightly closed ideal in $R$, we have
$(I+P) \cap (I+Q) = I+(P \cap Q)$ as ideals of $S$, by Lemma \ref{intersect}. 
Any element of the
maximal ideal of $S$, up to multiplication by units, is either in $Q$, or is
of the form $Y^h +q$, where $q \in Q$. Since $I$ cannot be contained in
$Q$, one of its generators has the form $Y^h +q$. Choosing the generator 
amongst these which has the least such positive value of $h$, and subtracting 
suitable multiples of this generator, we may assume that the other generators
are in $Q$. We then have $I=(Y^h+q_1,q_2,\dots,q_d)S$, where $q_i \in Q$, 
$h >0$, and $d=\dim R=\dim S/P$. By a similar argument we may write $P$ as 
$P=(Y^t+r_1,r_2, \dots, r_k)S$, where $r_i\in Q$. Since we are assuming
that the image of $I$ is a tightly closed ideal in $R$, Theorem \ref{contain}
shows that $I$ is not contained in  $P+Q$, and so we conclude $h < t$.

\medskip

We then have $Y^t+Y^{t-h}q_1 \in I \cap (P+Q)$, and so by Lemma
\ref{distribute}
$$
Y^t+Y^{t-h}q_1 \in (I \cap P)+(I \cap Q).
$$
By Lemma \ref{regseq}, $I \cap P= IP$ and consequently
$$
Y^t \in IP +Q = (Y^t+r_1)(Y^h+q_1)+Q=Y^{t+h}+Q.
$$ 
However this is impossible since $h > 0$. 
\end{proof}

\begin{rem}                                             
Note that in the proof above we used that if $\overline I$ is a tightly closed 
ideal, we must have $I + (P \cap Q) = (I+P)\cap (I+Q)$, and then showed that 
this cannot hold when $\overline I$ is generated by a system of parameters for 
$R$ in the case $S/P$ is Cohen--Macaulay, and $S/Q$ is a discrete valuation 
ring. When $S/P$ is not Cohen--Macaulay, this approach no longer works as seen
from the following example. 
\end{rem}

\begin{ex}
Let $S=K[[T,X,Y,Z]]$, and consider the two prime ideals $Q=(T,X,Y)$ and 
$P=(TY-XZ,T^2X-Z^2,TX^2-YZ,X^3-Y^2)$. Then $S/Q$ is a discrete 
valuation ring, although $S/P$ is not Cohen--Macaulay. To see this, observe 
that 
$$
S/P \cong K[[U^2, U^3, UT, T]] \subseteq K[[T, U]]
$$ 
where $T$ and $U$ are indeterminates and $x \mapsto U^2$,
$y \mapsto U^3$, $z \mapsto UT$ and $t \mapsto T$. (Lower case letters denote
the images of the corresponding variables.)

\medskip

Then $R=S/(P \cap Q)$ is a \NE ring and the image of $I=(Z,X-T)$ in $R$ is 
$\overline I=(z,x-t)$ which is generated by a system of parameters
for $R$. We shall see that $I+(P \cap Q) = (I+P)\cap (I+Q)$.

\medskip

Since $I+Q = (T,X,Y,Z)$ is just the maximal ideal of $S$, we get that 
$(I+P)\cap (I+Q) = I+P = (Z, X-T, XY, X^3, Y^2)$.
It can be verified (using Macaulay, or even otherwise) that 
$$
P \cap Q = (TY-XZ, TX^2-YZ, X^3-Y^2, T^3X-TZ^2)
$$ 
and so 
$$
I + (P \cap Q) = (Z, X-T, XY, X^3, Y^2) = (I+P)\cap (I+Q).
$$ 

\medskip

For the ring $R$, although it does not follow from any of the earlier results,
we can show that no system of parameters generates a tightly closed ideal.
\label{badintersect}
\end{ex}

We can actually prove the graded analogue of Theorem \ref{CMDVR} without the
requirement that $S/P$ is Cohen--Macaulay.

\begin{thm}
Let $S=K[X_1, \dots, X_n, Y]$ with ideal $Q=(X_1,\dots, X_n)S$, and 
$P$ a homogeneous unmixed ideal with $\dim S/P \ge 2$. Then no homogeneous
system of parameters of the ring $R=S/(P \cap Q)$ generates a tightly closed
ideal.
\label{NCMDVR}
\end{thm}

\begin{proof} 
Let $I$ be an ideal of $S$ generated by homogeneous elements which map to a 
system of parameters in $R$, and assume that the image of $I$ is a tightly
closed ideal of $R$.  

\medskip

As in the proof of Theorem \ref{CMDVR}, there is no loss of generality in
taking as homogeneous generators for $I$, the elements 
$Y^h+q_1,\ q_2, \dots, q_d$ where $q_i \in Q$, and $h >0$, and for $P$ the 
elements $Y^t+r_1,\ r_2, \dots, r_k$, where $r_i\in Q$. One can easily formulate 
a graded analogue of Theorem \ref{contain} and then since we are assuming that
the image of $I$ is a tightly closed ideal in $R$, it follows that $I$ is not 
contained in  $P+Q$. Hence we conclude $h < t$.

The assumption implies that 
$$
Y^t+r_1 \in (I+P) \cap (I+Q) = I+(P \cap Q) = I+(r_2, \dots, r_k) + (Y^t+r_1)Q
$$
and so $Y^t+r_1 \in I+(r_2, \dots, r_k)$. Hence 
$$
I+(P \cap Q) = I+P = I+(r_2, \dots, r_k).
$$

\medskip

If $S/P$ is Cohen--Macaulay, the proof is identical to that of Theorem
\ref{CMDVR}, and so we may assume $S/P$ is not Cohen--Macaulay. Consequently
$(IS/P)^*$ is strictly bigger that $IS/P$. Let $F \in S$ be a homogeneous element
such that its image is in $(IS/P)^*$ but not in $IS/P$. Note that if $F \in I+Q$,
then $\overline F \in (IR)^*$, and so $F \in I+P$, a contradiction. Hence
we conclude that $F \notin I+Q = (Y^h, X_1, \dots, X_n)$ and so $F = Y^i + G$ 
where $i < h$ and $G \in Q$.

\medskip

Next note that 
$Y^{h-i}F = Y^h + GY^{h-i} \in I+(P \cap Q) = I+(r_2, \dots, r_k)$, and so 
$GY^{h-i} - q_1 \in (q_2, \dots, q_d, r_2, \dots, r_k)$. We then have
$$
\overline F  \in (IS/P)^* = ((Y^h+GY^{h-i}, q_2, \dots, q_d)S/P)^*
 = (Y^{h-i}F, q_2, \dots, q_d)S/P)^*. 
 $$
By a degree argument, we see that 
$\overline F  \in ((q_2, \dots, q_d)S/P)^*$. However this means that 
$\overline F$ is in the radical of the ideal $(q_2, \dots, q_d)S/P$, which
contradicts the fact that $(FY^{h-i}, q_2, \dots, q_d) = IS/P$ is primary to
the homogeneous maximal ideal of $S/P$. 
\end{proof}

\section{NE closure}

For Noetherian rings of characteristic $p$ we shall define a new closure
operation on ideals, the {\it NE closure}, which will agree with tight closure 
when the ring is equidimensional. In \NE local rings, tight closure no longer 
has the so--called {\it colon--capturing} property, and the main point of NE 
closure is to recover this property. This often forces the NE closure of an 
ideal to be larger than its tight closure and at times even larger than its 
radical, see Example \ref{badcolon}. More precisely let $(R,m)$ be an 
excellent local ring with a system of parameters $x_1,\dots,x_n$. Then when $R$ is 
equidimensional we have $(x_1,\dots,x_k):x_{k+1} \subseteq (x_1,\dots,x_k)^*$, 
but this does not hold in general. The NE closure (denoted by
$I\bs$ for an ideal $I \subseteq R$) will have the property that
$(x_1,\dots,x_k):x_{k+1} \subseteq (x_1,\dots,x_k)\bs$.

\begin{defn}
We shall say that a minimal prime ideal $P$ of a ring $R$ is {\it absolutely
minimal} if $\dim R/P = \dim R$. When Spec $R$ is connected, $R\dot$ shall
denote the complement in $R$ of the union of all the absolutely minimal
primes. If $R= \prod R_i$, we define $R\dot = \prod R_i\dot$.
The {\it NE closure} $I\bs$ of an ideal $I$ is given by 

$I\bs=\{x \in R :$ there exists 
$c \in R\dot$ with $cx^{[q]} \in I^{[q]}$ for all  $q \gg 0$\}.
\end{defn}

\medskip

The following proposition and its proof are analogous to the 
statements for tight closure in equidimensional rings, see 
\cite[Theorem 4.3]{HHbasec}.
 
\begin{prop}
Let $R$ be a complete local ring of characteristic $p$, with a system of 
parameters $x_1,\dots,x_n$. Then 

\item $(1)$\quad $(x_1,\dots,x_k):x_{k+1} \subseteq (x_1,\dots,x_k)\bs$.

\item $(2)$\quad $(x_1,\dots,x_k)\bs:x_{k+1} = (x_1,\dots,x_k)\bs$.

\item $(3)\quad $If $(x_1,\dots,x_{k+1})\bs=(x_1,\dots,x_{k+1})$, then 
$(x_1,\dots,x_k)\bs=(x_1,\dots,x_k)$.

\item $(4)\quad $If $(x_1,\dots,x_n)\bs=(x_1,\dots,x_n)$ or 
$(x_1,\dots,x_{n-1})\bs=(x_1,\dots,x_{n-1})$, then $R$ is Cohen--Macaulay.
\label{necapturing}
\end{prop}

\begin{proof} 
(1)\quad We may represent $R$ as a module-finite extension of a regular subring $A$ 
of the form $A=K[[x_1,\dots,x_n]]$ where $K$ is a field. Let $t$ be the
torsion free rank of $R$ as an $A$--module, and consider $A^t \subseteq R$.
Then $R/A^t$ is a torsion $A$--module and there exists $c \in A$, nonzero, such 
that $cR \subseteq A^t \subseteq R$. We note that $c$ cannot be in any
absolutely minimal prime $P$ of $R$, since for any such $P$, $R/P$ is of
dimension n and is module-finite over $A/{A \cap P}$, and so $A \cap P=0$.
Now if $u \in (x_1,\dots,x_k):x_{k+1}$ then for some $r_i \in R$, 
$ux_{k+1}=\sum_{i=1}^k r_ix_i$. Taking $q$th powers, and multiplying by $c$ we
get $cu^q x_{k+1}^q=\sum_{i=1}^k cr_i^qx_i^q$. But now $cu^q$ and each of
$cr_i^q$ are in $A^t$ and $x_i^q$ form a regular sequence on $A^t$. 
Hence $cu^q \in (x_1^q,\dots,x_k^q)$ and so $u \in (x_1,\dots,x_k)\bs$.

(2)\quad If $ux_{k+1} \in (x_1,\dots,x_k)\bs$ then for some $c_0 \in R\dot$, 
$c_0(ux_{k+1})^q \in (x_1^q,\dots,x_k^q)$ for all  sufficiently large $q$, 
i.e., $c_0u^q x_{k+1}^q=\sum_{i=1}^kr_ix_i^q$ for $q \gg 0$. 
Multiplying this by our earlier choice of $c$, we again have a relation on
$x_i^q$'s with coefficients in $A^t$, namely 
$cc_0u^q x_{k+1}^q=\sum_{i=1}^k cr_ix_i^q$ for $q \gg 0$, and so 
$cc_0u^q \in (x_1^q,\dots,x_k^q)$ for $q \gg 0$. Since $cc_0 \in R\dot$ we get 
$u \in (x_1,\dots,x_k)\bs$.

(3)\quad Let $J=(x_1,\dots,x_k)$. Then $J\bs \subseteq (x_1,\dots,x_{k+1})$
and so $J\bs \subseteq J+x_{k+1}R$. If $u \in J\bs$, $u=j+x_{k+1}r$
for $j \in J$ and $r \in R$. This means $r \in J\bs :x_{k+1}$ which equals
$J\bs$ by (2). Hence we get $J\bs =J+x_{k+1}J\bs$. Now by Nakayama's lemma we
get $J\bs=J$.

(4)\quad This follows from (2) and (3) since, under either of the hypotheses, the
system of parameters $x_1,\dots,x_n$ is a regular sequence.
\end{proof}

The above proposition, coupled with results on F--rationality, gives us the
following theorem:

\begin{thm}
Let $R$ be a complete local ring of characteristic $p$, with a system of 
parameters which generates a NE--closed ideal. Then $R$ is F--rational.
\label{neclosed}
\end{thm}

\begin{proof} 
From the previous proposition the ring is Cohen--Macaulay, and in particular,
equidimensional. For equidimensional rings, tight closure agrees with NE
closure, and the result follows from Theorem \ref{longlist} (2).
\end{proof}

We shall extend this result to excellent local rings once we develop the
theory of test elements for NE closure. The following proposition
lists some properties of NE closure.

\begin{prop}
Let $R$ be a ring of characteristic $p$, and $I$ an ideal of $R$.

\item $(1)$\quad $0\bs$ is the intersection of the absolutely minimal prime ideals of $R$.

\item $(2)$\quad If $I=I\bs$ then for any ideal $J$, $(I:J)\bs=I:J$.

\item $(3)$\quad If $R= \prod R_i$ and $I= \prod I_i$, then $I\bs=\prod I_i\bs$.

\item $(4)$\quad For rings $R$ and $S$ and a homomorphism $h:R \to S$ satisfying 
the condition $h(R\dot) \subseteq S\dot$, we have $h(I\bs) \subseteq (IS)\bs$.

\item $(5)$\quad $x \in I\bs$ if and only if $\overline x \in (IR/P)\bs$ for every
absolutely minimal prime ideal $P$ of $R$.
\label{nebasics}
\end{prop}

\begin{proof}  (1), (2), (3) and (4) follow easily from the definitions. For (5)
note that if $P$ is absolutely minimal, $h:R \to R/P$ meets the condition of 
(4), so 
$x \in I\bs$ implies that its image is in the NE closure of $IR/P$. For the
converse, fix for every absolutely minimal $P_i$, $d_i \notin P_i$ but in every
other minimal prime of $R$. If $\overline x \in (IR/P_i)\bs$ for every
absolutely minimal $P_i$, 
then there exist elements $\overline c_i$ with 
$\overline {c_ix^q} \in (IR/P_i)^{[q]}$. We can lift each $\overline c_i$ to 
$c_i \in R$ with $c_i \notin P_i$. Then $c_ix^q \in I^{[q]}+P_i$ for all $i$,
for sufficiently large $q$.
Multiplying each of these equations with the corresponding $d_i$, we get
$c_id_ix^q \in I^{[q]}+{\mathfrak{N}}$, since $d_iP_i$ is a subset of every minimal
prime ideal and so is in the nilradical, ${\mathfrak{N}}$. If ${\mathfrak{N}}^{[q']}=0$,
taking $q'$ powers of these equations gives us $(c_id_i)^{q'}x^q \in I^{[q]}$ 
for all $i$, for sufficiently large $q$. Set $c=\sum (c_id_i)^{q'}$. By our 
choice of $c_i$'s and $d_i$'s, $c\in R\dot$, and the above equations put
together give us $cx^q \in I^{[q]}$ for all sufficiently large $q$. 
\end{proof}

We note that NE closure need not be preserved once we localize, i.e., it is
quite possible that $x \in I\bs$, but $x \notin (IR_P)\bs$. Examples of this 
abound in  non--equidimensional rings, but there are some positive 
results about NE closure being preserved under certain maps which we examine 
in the next few propositions.

\begin{prop}
If $h:(R,m) \to (S,n)$ is a faithfully flat homomorphism of local rings then
for an ideal $I$ of $R$, if $x\in I\bs$, then its image $h(x)$ is in $(IS)\bs$. 
In particular if $\HR$ denotes the completion of $R$ at its maximal ideal, 
$x \in I\bs$ implies $x \in (I\HR)\bs$.
\end{prop}

\begin{proof} 
By Proposition \ref{nebasics} (4), it suffices to check that 
$h(R\dot) \subseteq S\dot$.
This is equivalent to the assertion that the contraction of every absolutely
minimal prime of $S$ is an absolutely minimal prime of $R$. Now let $P$ be an
absolutely minimal prime of $S$, and $p$ denote its contraction to $R$.
Then since $R \to S$ is faithfully flat, by a change of base, so is $R/p \to
S/pS$. This gives dim $S/pS=$ dim $R/p$ $+$ dim $S/mS$. Also, faithful flatness
of $h$ implies that dim $S=$ dim $R$ $+$ dim $S/mS$. But $P$ was an absolutely 
minimal prime of $S$, so dim $S=$ dim $S/P=$ dim $S/pS$, since $pS \subseteq
P$. Putting these equations together, we get dim $R/p=$ dim $R$, and so $p$ is
an absolutely minimal prime of $R$.
\end{proof}

\begin{prop}
Let $R$ and $S$ be Noetherian rings of characteristic $p$, and $R \to S$ a 
homomorphism such that for every absolutely minimal prime $Q$ of $S$, its 
contraction to $R$, $Q^c$, contains an absolutely minimal prime of $R$. Assume
one of the following holds:

\item $(1)$\quad $R$ is finitely generated over an excellent local ring, or is 
F--finite, or

\item $(2)$\quad $R$ is locally excellent and $S$ has a locally stable test 
element, (or $S$ is local), or

\item $(3)$\quad $S$ has a completely stable test element (or $S$ is a complete 
local ring).

Then if $x \in I\bs$ for $I$ an ideal of $R$, the image of $x$ in $S$ is in 
$(IS)\bs$.
\end{prop}

\begin{proof} 
It suffices to check $x \in (IS/Q)\bs$ for every absolutely minimal primes $Q$
of $S$, by Proposition \ref{nebasics} (5).  But $(IS/Q)\bs=(IS/Q)^*$ since 
$S/Q$ is equidimensional. If $P \subseteq Q^c$ is an absolutely minimal prime of 
$R$, then $x \in I\bs$ implies $\overline x \in (IR/P)\bs=(IR/P)^*$. The result 
now follows by applying \cite[Theorem 6.24]{HHbasec} to the map $R/P \to S/Q$.
\end{proof}

\section{NE--test elements}

\begin{defn}
We shall say $c\in R\dot$ is a {\it $q'$--weak NE--test element for $R$} if 
for all ideals $I$ of $R$ and $x\in I\bs$, $cx^q\in I^{[q]}$ for all $q \ge
q'$. We may often use the phrase {\it weak NE--test element} and suppress the
$q'$.

\medskip

For a local ring $(R,m)$, $c\in R\dot$ is a 
{\it weak completion stable NE--test element for $R$} if it is a weak 
NE--test element for $\HR$, the completion of $R$ at its maximal ideal. 
\end{defn}

Our definition of a completion stable weak NE--test element is different
from the notion of a completely stable weak test element for tight closure, 
where it is required
that the element serve as a weak test element in the completion of every local
ring of $R$. The reason for this, of course, is that localization is no longer
freely available to us, since $R\dot$ often does not map into $(R_P)\dot$.

\medskip

Note also that since $\HR$ is faithfully flat over $R$, a weak completion 
stable  NE--test element for $R$ is also a weak NE--test element for $R$.

\begin{prop}
If for every absolutely minimal prime $P$ of $R$, $R/P$ has a weak test
element, then $R$ has a weak NE--test element. 
\end{prop}

\begin{proof} 
Fix for every absolutely minimal prime $P_i$ an element $d_i$ not in $P_i$ but 
in every other 
minimal prime of $R$. Let ${\mathfrak{N}}$ denote the nilradical of $R$ and fix $q'$
such that ${\mathfrak{N}}^{[q']}=0$. If $\overline c_i$ is a weak test element for 
$R/P_i$, we may pick $c_i \notin P_i$ which maps to it under $R \to R/P$.
We claim $c=\sum (c_id_i)^{q'}$ is a weak NE--test element for $R$. If 
$x \in I\bs$, then $\overline x \in (IR/P_i)\bs$ for all $P_i$ absolutely
minimal. Since $\overline c_i$ is a weak test element for $R/P_i$, we have
$\overline {c_ix^q} \in (IR/P_i)^{[q]}$ for all $i$, for sufficiently large
$q$, i.e., $c_ix^q \in I^{[q]} + P_i$. Multiplying this by $d_i$, summing over
all $i$ and taking the $q'$ power as in the proof of Proposition
\ref{nebasics} (5), we
get that $cx^q \in I^{[q]}$. It is easy to see that $c\in R\dot$ and so is a 
weak NE--test element. 
\end{proof}

\begin{prop}
Every excellent local ring of characteristic $p$ has a weak completion stable 
NE--test element.
\end{prop}

\begin{proof} 
If $R$ is an excellent local domain, it has a completely stable weak test
element, see \cite[Theorem 6.1]{HHbasec}. Hence each $R/P_i$ for $P_i$ 
absolutely minimal, has a completely 
stable weak test element, say  $\overline c_i$. Let $c_i$, $d_i$, $q'$ and $c$
be as in the proof of the previous proposition. If $x \in (I\HR)\bs$, we have
$\overline x \in (I\HR/P_i\HR)\bs$. Since $R/P_i$ is equidimensional and
excellent, its completion $\HR/P_i\HR$ is also equidimensional. (We use here
the fact that the completion of a universally catenary equidimensional local 
ring is again equidimensional, \cite[Page 142]{HIO}). Hence NE closure agrees
with tight closure in $\HR/P_i\HR$, and we get $\overline x \in 
(I\HR/P_i\HR)^*$. This gives $\overline {c_ix^q} \in (I\HR/P_i\HR)^{[q]}$
for all $i$, for sufficiently large $q$. As in the previous proof, we then get
that $cx^q \in (I\HR)^{[q]}$, and so is a weak completion stable NE--test 
element.
\end{proof}

We can now extend Theorem \ref{neclosed} to the case where $R$ is 
excellent local, without requiring it to be complete.

\begin{thm}
Let $(R,m)$ be an excellent local ring of characteristic $p$ with a system of 
parameters $x_1,\dots,x_n$. Then  
\item $(1)$\quad $(x_1,\dots,x_k):x_{k+1} \subseteq (x_1,\dots,x_k)\bs$.

\item $(2)$\quad $(x_1,\dots,x_k)\bs:x_{k+1} = (x_1,\dots,x_k)\bs$.

\item $(3)$\quad If $(x_1,\dots,x_{k+1})\bs=(x_1,\dots,x_{k+1})$, then 
$(x_1,\dots,x_k)\bs=(x_1,\dots,x_k)$.

\item $(4)$\quad If $(x_1,\dots,x_{n-1})\bs=(x_1,\dots,x_{n-1})$ then $R$ is 
Cohen--Macaulay. 

\item $(5)$\quad If $(x_1,\dots,x_n)\bs=(x_1,\dots,x_n)$ then $R$ is F--rational.
\end{thm}  

\begin{proof} 
Since $R$ has a weak completion stable NE--test element, if there is a
counterexample to any of the above claims, we can preserve this while mapping
to $\HR$. But all of the above are true for complete local rings as follows
from Proposition \ref{necapturing} and Theorem \ref{neclosed}.
\end{proof}

\begin{ex}
Let $R=K[[X,Y,Z]]/(X)\cap(Y,Z)$. Then $y,x-z$ is a system of parameters for $R$
and $0:_R(y) =(x)$. That tight closure fails here to \lq\lq capture
colons\rq\rq \ is seen from the fact that $x \notin 0^*=0$. However 
$0\bs =  (x)$, and we certainly have $0:_R(y) \subseteq 0\bs$.
\label{badcolon}
\end{ex}

\medskip

\section*{Acknowledgments}

It is a pleasure to thank Melvin Hochster for many valuable discussions.

\end{document}